\documentclass[a4paper]{amsart}

\usepackage{amsthm, amsfonts, amssymb, amsmath, latexsym, enumerate, times}
\usepackage[latin1]{inputenc}
\usepackage{mathrsfs}
\usepackage{mathtools}
\usepackage{letltxmacro}
\LetLtxMacro\orgvdots\vdots
\LetLtxMacro\orgddots\dots

\makeatletter
\DeclareRobustCommand\vdots{%
	\mathpalette\@vdots{}%
}
\newcommand*{\@vdots}[2]{%
	\sbox0{$#1\cdotp\cdotp\cdotp\m@th$}%
	\sbox2{$#1.\m@th$}%
	\vbox{%
		\dimen@=\wd0 %
		\advance\dimen@ -3\ht2 %
		\kern.5\dimen@
		\dimen@=\wd2 %
		\advance\dimen@ -\ht2 %
		\dimen2=\wd0 %
		\advance\dimen2 -\dimen@
		\vbox to \dimen2{%
			\offinterlineskip
			\copy2 \vfill\copy2 \vfill\copy2 %
		}%
	}%
}
\DeclareRobustCommand\ddots{%
	\mathinner{%
		\mathpalette\@ddots{}%
		\mkern\thinmuskip
	}%
}
\newcommand*{\@ddots}[2]{%
	\sbox0{$#1\cdotp\cdotp\cdotp\m@th$}%
	\sbox2{$#1.\m@th$}%
	\vbox{%
		\dimen@=\wd0 %
		\advance\dimen@ -3\ht2 %
		\kern.5\dimen@
		\dimen@=\wd2 %
		\advance\dimen@ -\ht2 %
		\dimen2=\wd0 %
		\advance\dimen2 -\dimen@
		\vbox to \dimen2{%
			\offinterlineskip
			\hbox{$#1\mathpunct{.}\m@th$}%
			\vfill
			\hbox{$#1\mathpunct{\kern\wd2}\mathpunct{.}\m@th$}%
			\vfill
			\hbox{$#1\mathpunct{\kern\wd2}\mathpunct{\kern\wd2}\mathpunct{.}\m@th$}%
		}%
	}%
}
\makeatother

\usepackage{array}
\usepackage{comment}
\usepackage{paralist}
\usepackage{xcolor}

\newtheorem{theorem}{Theorem}
\newtheorem{lemma}[theorem]{Lemma}

\newtheorem{proposition}[theorem]{Proposition}

\newtheorem{question}[theorem]{Question}

\theoremstyle{definition}

\newcommand{\bbQ}{{\mathbb Q}}
\newcommand{\bbC}{{\mathbb C}}

\def\le{\leqslant}
\def\ge{\geqslant}

\usepackage{eurosym}

\begin{document}
 
\title{A question about generalized nonvanishing}
 
\author{Claudio Fontanari}

\address{Claudio Fontanari, Dipartimento di Matematica, Universit\`a degli Studi di Trento, Via Sommarive 14, 38123 Povo, Trento}
\email{claudio.fontanari@unitn.it}

\subjclass[2010]{14E30} 
 
\keywords{Nonvanishing, Generalized Nonvanishing, Effective Cone, Pseudoeffective Cone.}
 
\centerline{}
\begin{abstract} 
Let $(X,\Delta)$ be a projective, $\bbQ$-factorial log canonical pair and let $L$ be a pseudoeffective $\bbQ$-divisor on $X$ such that 
$K_X + \Delta + L$ is pseudoeffective. Is there an effective $\bbQ$-divisor $M$ on $X$ such that $K_X + \Delta + L$ is numerically equivalent to $M$? 
We are not aware of any counterexamples, but the answer is not completely clear even in the case of surfaces. 
\end{abstract}

\maketitle 

\section{Introduction}

In 2014, Birkar and Hu asked the following question in the context of the theory of pairs polarized by a nef divisor (see \cite{BH}, Question 3.5): 

\begin{question}\label{nef}
Let $(X,\Delta)$ be a projective, $\bbQ$-factorial, log canonical pair and let $L$ be a nef $\bbQ$-divisor on $X$ such that 
$K_X + \Delta + L$ is pseudoeffective. Is there an effective $\bbQ$-divisor $M$ on $X$ such that 
$K_X + \Delta + L$ is numerically equivalent to $M$?
\end{question}

The formulation in terms of numerical (and not linear) equivalence turns out to be necessary: indeed, a non-torsion numerically trivial divisor 
on an elliptic curve is nef but has negative Kodaira dimension. In spite of the skepticism shown by Birkar and Hu (see \cite{BH}, p. 212: 
\emph{Most probably, the answer is no. However, there are interesting cases in which the answer is yes}), there has been recently quite remarkable
progress towards an affirmative answer. In particular, Lazi\'c and Peternell have established that it would follow from the standard conjectures of the 
Minimal Model Program (namely, termination of flips, Abundance Conjecture and Semiampleness Conjecture), at least in the case $(X,\Delta)$ klt and 
$K_X + \Delta$ pseudoeffective (see \cite{LP}, Theorem A). Furthermore, Question \ref{nef} has affirmative answers also for surfaces 
(see \cite{HL}, Theorem 1.5) and for threefolds with nef anticanonical divisor (see \cite{L}, Theorem A). 

On the other hand, if one does not look at Question \ref{nef} as a statement on $X$ polarized by the divisor $L$, but rather as a generalization of the classical nonvanishing for the canonical divisor of $X$, one may wonder whether the assumption that $L$ is 
nef is really essential. Indeed, here we pose the following stronger question, where nefness is replaced by pseudoeffectivity: 

\begin{question}\label{psef}
Let $(X,\Delta)$ be a projective, $\bbQ$-factorial log canonical pair and let $L$ be a pseudoeffective $\bbQ$-divisor on $X$ such that 
$K_X + \Delta + L$ is pseudoeffective. Is there an effective $\bbQ$-divisor $M$ on $X$ such that $K_X + \Delta + L$ is numerically equivalent to $M$?
\end{question}

In particular, an affirmative answer would imply that on a variety $X$ such that $-K_X$ is pseudoeffective (for instance, a K-trivial variety) 
with $h^1(X, \mathcal{O}_X)=0$, so that linear and numerical equivalence coincide, any pseudoeffective $\bbQ$-divisor is indeed effective. 
To the best of our knowledge, there is not an explicit expectation in this sense but there are no counterexamples as well. In dimension $2$, 
it is known that if the anticanonical divisor of a regular surface is nef and non-torsion then the effective cone is closed (see \cite{Bo}, 
Proposition 6.2); the same conclusion holds also in the case of $\mathbb{P}^3$ blown up in eight very general points, where the anticanonical 
divisor is nef but not semiample (see \cite{SX}, Theorem 1.4).  
 
We are going to present the following evidence in dimension $2$ towards an affirmative answer to Question \ref{psef}. 

\begin{proposition}\label{pseudoeffective}
Let $X$ be a smooth projective surface and let $L$ be a pseudoeffective divisor on $X$. 
Assume that $K_X$ is pseudoeffective. Then $K_X + L$ is numerically equivalent to an effective $\bbQ$-divisor.
\end{proposition}

\begin{theorem}\label{regular}
Let $X$ be a smooth projective surface and let $L$ be a pseudoeffective divisor on $X$ such that $K_X + L$ is pseudoeffective. 
Assume that $h^1(X, \mathcal{O}_X)=0$. Then a multiple of $K_X + L$ is effective.
\end{theorem}

If $K_X$ is not pseudoeffective, i.e. if $X$ is uniruled, or the irregularity $q = q(X) = h^1(X, \mathcal{O}_X)$ is strictly 
positive, then the answer seems to be more elusive. An old result by Sakai (see \cite{S}, Proposition 8) applies to an arbitrary surface, 
but requires that the pseudoeffective divisor $L$ is a reduced effective divisor.

\begin{proposition}\label{reduced} \emph{(Sakai)}
Let $X$ be a smooth projective surface and let $D$ be a reduced effective divisor on $X$ such that $K_X + D$ is pseudoeffective. 
Then a multiple of $K_X + D$ is effective.
\end{proposition}

We are able to address a broader class of effective divisors under a few additional technical assumptions. Recall that an effective 
divisor $D$ on $X$ is $0$-connected if for every decomposition $D=D_1+D_2$ into two effective divisors we have $D_1.D_2 \ge 0$. 

\begin{theorem}\label{connected}
Let $X$ be a smooth projective ruled surface of irregularity $q \ge 2$ and let $D$ be a $0$-connected effective divisor on $X$ such that $K_X + D$ is nef. 
Then a multiple of $K_X + D$ is effective.
\end{theorem}

We work over the complex field $\bbC$.

\medskip

{\bf Acknowledgements:} 
The author is grateful to Marco Andreatta, Edoardo Ballico, Paolo Cascini, Ciro Ciliberto, Vladimir Lazi\'c, Roberto Pignatelli 
and Roberto Svaldi for helpful conversations on these and related topics. 
The author is a member of GNSAGA of INdAM (Italy). 

\section{The proofs}

\noindent
\textit{Proof of Proposition \ref{pseudoeffective}.}  
By the Zariski decomposition (see for instance \cite{B}, Theorem 14.14) we have $L = P(L) + N(L)$, 
with $P(L)$ nef and $N(L)$ effective. In particular, $K_X + P(L)$ is pseudoeffective and by \cite{HL}, Theorem 1.5, $K_X + P(L)$ 
is numerically equivalent to an effective $\bbQ$-divisor. It follows that the same is true for $K_X + L = K_X + P(L) + N(L)$ as well. 
  
\qed

\noindent
\textit{Proof of Theorem \ref{regular}.} By the Zariski decomposition, now applied to the pseudoeffective divisor $K_X + L$,
we have $K_X + L = P + N$, with $P = P(K_X+L)$ nef, $N = N(K_X+L)$ effective and $P.N = 0$.
If $P^2 > 0$ then $P$ is big and a suitable multiple of $K_X+L$ is effective. If $P$ is numerically equivalent to zero, then $K_X+L$ is 
numerically (and also linearly, since $h^1(X, \mathcal{O}_X)=0$) equivalent to the effective divisor $N$. 
Hence we may assume both $P^2 = 0$ and $P$ not numerically equivalent to zero. In 
particular, $-P$ is not pseudoeffective and $h^2(X, mP) = h^0(X, K_X - mP) = 0$ for $m >>0$ (see for instance \cite{B}, Lemma 14.6). 
As a consequence, the Riemann-Roch formula takes the simple form:
$$
h^0(X, mP) = h^1(X, mP) + \chi(\mathcal{O}_X) - \frac{1}{2} mP.K_X \ge 1 - \frac{1}{2} mP.K_X  
$$  
by the assumption $h^1(X, \mathcal{O}_X)=0$. Next, we claim that $P.K_X \le 0$. Indeed, we have $P.(K_X+L) = P.(P+N) = P^2 + P.N = 0$ 
and $P.L \ge 0$ since $P$ is nef and $L$ is pseudoeffective. Hence we conclude that $h^0(X, m(K_X + L)) \ge h^0(X, mP) > 0$. 

\qed 

\noindent

\begin{lemma}\label{ruled}
Let $X$ be a smooth projective ruled surface and let $L$ be a pseudoeffective divisor on $X$ such that $K_X + L$ is nef. 
If $q \ge 1$ assume also that $(K_X+L).L > 0$. Then a multiple of $K_X + L$ is effective.
\end{lemma}

\proof If we apply Riemann-Roch to the divisor $m(K_X +L)$ , where $m$ is a positive integer, we obtain: 
\begin{eqnarray*}
h^0(X, m(K_X +L)) - h^1(X, m(K_X +L)) + h^2(X, m(K_X +L)) \\
= \chi(\mathcal{O}_X) + \frac{1}{2} \left( m(m-1)(K_X+L)^2 + m(K_X+L).L \right)
\end{eqnarray*}
and by Serre duality we have $h^2(X, m(K_X +L))= h^0(X, - L - (m-1)(K+L)) =0$ since $L + (m-1)(K+L)$ is a nonzero pseudoeffective divisor. 
Moreover, we have $m(m-1)(K_X+L)^2 \ge 0$ and $m(K_X+L).L \ge 0$ since $K_X+L$ is assumed to be nef. It follows that 
$h^0(X, m(K_X +L)) > 0$ as soon as either $\chi(\mathcal{O}_X) = 1-q > 0$ or $(K_X+L).L > 0$. 

\qed

\noindent
\textit{Proof of Theorem \ref{connected}.} In order to apply Lemma \ref{ruled} we need to check that $(K_X+D).D > 0$. Assume by contradiction 
that $(K_X+D).D = 0$. By adjunction, the arithmetic genus of $D$ is $p_a(D)=1$. We claim that every irreducible component $C$ in the support of 
$D$ satisfies $p_a(C) \le p_a(D)$. Indeed, we have $2p_a(D)-2 = (K_X+D).D = (K_X+D).(C+D-C)= (K_X+D).C + (K_X+D).(D-C) \ge (K_X+C).C + 
(D-C).C \ge (K_X+C).C = 2p_a(C)-2$, where the two inequalities hold since $K_X+D$ is nef and $D$ is $0$-connected. On the other hand, we 
claim that at least one irreducible component $C$ in the support of $D$ is not contained in a fiber $F$ of the projection of $X$ onto its base curve $B$. 
Indeed, otherwise we would have $D.F=0$, contradicting the nefness of $K_X+D$ since $K_X.F = -2$. Summing up, there exists a curve $C$ on $X$ 
with $p_a(C) \le 1$ not contained in a fiber, hence projecting onto the base curve $B$. If follows that $q(X) = p_a(B) \le p_a(C) \le 1$, contradicting the assumption $q \ge 2$. 

\qed

\end{document}